\documentclass{article}
\usepackage{amsmath,amsthm,amssymb}

\def\Ex{{\mathbb E}}
\def\Pr{{\mathbb P}}
\def\ind{\mathbf I}
\def\er{{\mathbb R}}
\def\ve{\varepsilon}

\def\Cov{\mathrm{Cov}}
\def\Var{\mathrm{Var}}
\def\Med{\mathrm{Med}}
\def\trip{|\!|\!|}

\newtheorem{thm}{Theorem}

\newtheorem{prop}[thm]{Proposition}

\newtheorem{cor}[thm]{Corollary}

\title{Moments of unconditional logarithmically concave vectors
\thanks{
Research partially supported by MNiSW Grant no. N N201 397437 and the Foundation for Polish 
Science.}}
\author{Rafa{\l} Lata{\l}a \thanks{Institute of Mathematics, University of Warsaw, 
Banacha 2, 02-097 Warszawa, Poland 
and Institute of Mathematics, Polish Academy of Sciences, 
ul. \'{S}niadeckich 8, 00-956 Warszawa, Poland, e-mail: rlatala@mimuw.edu.pl.}}
\date{}

\begin{document}

\maketitle

\begin{abstract}
We derive two-sided bounds for moments of linear combinations of coordinates od
unconditional log-concave vectors. 
We also investigate how well moments of such combinations may be approximated by moments of Gaussian random variables.
\end{abstract}

\section{Introduction}

The aim of this paper is to study moments of
linear combinations of coordinates of unconditional, log-concave vectors
$X=(X_1,\ldots,X_n)$. A nondegerate random vector $X$ is 
\emph{log-concave} if it has a density of the form $g=e^{-h}$, where 
$h\colon \er\rightarrow (-\infty,\infty]$ is a convex function. 
We say that a random vector $X$  is \emph{unconditional} if the distribution of $(\eta_1X_1,\ldots,\eta_nX_n)$ is the same as $X$ for any choice
of signs $\eta_1,\ldots,\eta_n$. 

Typical example of uncoditional log-concave vector is a vector distributed uniformly in an unconditional
convex body $K$, i.e. such convex body that $(\pm x_1,\ldots,\pm x_n)\in K$ whenever 
$(x_1,\ldots,x_n)\in K$.  

A random vector $X$ is called {\em isotropic} if it has identity
covariance matrix, i.e.\ $\Cov(X_i,X_j)=\delta_{i,j}$. Notice that uncoditional vector $X$ is isotropic if and
only if its coordinates have variance one, in particular if $X$ is unconditional with nondenerate coordinates
then the vector $(X_1/\Var^{1/2}(X_1),\ldots,X_n/\Var^{1/2}(X_n))$ is isotropic and unconditional.

In \cite{GK} Gluskin and Kwapie\'n derived two-sided estimates for moments of $\sum_{i=1}^n a_iX_i$ if $X_i$
are independent, symmetric random variables with log-concave tails (coordinates of log-concave vector have
log-concave tails). In Section 2 we derive similar result for arbitrary unconditional
log-concave vectors $X$.  

In \cite{K} Klartag obtained powerful Berry-Essen type estimates for isotropic, unconditional, 
log-concave vectors $X$, showing in particular that if  
$\sum_{i} a_i^2=1$ and all $a_i$'s are small then the distribution of $S=\sum_{i=1}^n a_iX_i$
is close to the standard Gaussian distribution ${\cal N}(0,1)$. In Section 3 we investigate how well moments
of $S$ may be approximated by moments of ${\cal N}(0,1)$.

\medskip

{\bf Notation.}  By $\ve_1,\ve_2,\ldots$ we denote a Bernoulli sequence, i.e.\
a sequence of independent symmetric variables taking values $\pm 1$. We assume that
the sequence $(\ve_i)$ is independent of other random variables.

For a random variable $Y$ and $p>0$ we write $\|Y\|_p=(\Ex|Y|^p)^{1/p}$. For a sequence
$(a_i)$ and $1\leq q<\infty$, $\|a\|_q=(\sum_{i}|a_i|^q)^{1/p}$ and 
$\|a\|_{\infty}=\max_{i}|a_i|$. We set $B_q^n=\{a\in \er^n\colon \|a\|_q\leq 1\}$, $1\leq q\leq \infty$.
By $(a_i^*)_{1\leq i\leq n}$ we denote the nonincreasing rearrangement
of $(|a_i|)_{1\leq i\leq n}$.

We use letter $C$ (resp.\ $C(\alpha)$) for universal constants (resp.\ constants 
depending only on parameter $\alpha$). Value of a constant $C$ may differ at each
occurence. Whenever we want to fix the value of an absolute constant we will use letters 
$C_1,C_2,\ldots$. For two functions $f$ and $g$ we write $f\sim g$ to signify that 
$\frac{1}{C}f\leq g\leq Cf$.

\section{Estimation of moments}

It is well known and easy to show that if $X$ has a uniform distribution over a symmetric convex body
$K$ in $\er^n$ then for any $p\geq n$,
$\|\sum_{i\leq n} a_iX_i\|_p\sim \|a\|_{K^o}=\sup\{|\sum_{i\leq n} a_ix_i|\colon x\in K\}$. Our first
proposition generalizes this statement to arbitrary log-concave symmetric distributions.

\begin{prop}
\label{largep}
Suppose that $X$ has a symmetric $n$-dimensional log-concave distribution with the density $g$. 
Then for any $p\geq n$ we have
\[
\Big\|\sum_{i=1}^n a_iX_i\Big\|_p\sim \|a\|_{K_p^o},
\]
where
\[
K_p:=\{x\colon g(x)\geq e^{-p}g(0)\} \quad \mbox{and}\quad
\|a\|_{K_p^\circ}=\sup\Big\{\sum_{i=1}^n a_ix_i\colon x\in K_p\Big\}.
\]
\end{prop}

\begin{proof}
First notice that there exists an absolute constant $C_1$ such that
\[
\Pr(X\in C_1K_p)\geq 1-e^{-p}\geq \frac{1}{2}.
\]
For $n\leq p\leq 2n$ this follows by Corollary 2.4 and Lemma 2.2 in \cite{KM}.
For $p\geq 2n$ we may either adjust arguments from \cite{KM} or take any log-concave symmetric 
$m=\lfloor p\rfloor -n$ dimensional vector $Y$ independent of $X$ with density $g'$ and consider
the set $K'=\{(x,y)\in \er^n\times \er^m\colon g(x)g'(y)\geq e^{-p}g(0)g'(0)\}$. Then $K_p$
is a central $n$-dimensional section of $K'$, hence
$\Pr(X\in C_1K_p)\geq \Pr((X,Y)\in C_1\tilde{K'})\geq 1-e^{-p}$.

Observe that for any $z\in K_p$,
\[
\Big|\Big\{x\in K_p \colon \Big|\sum_{i=1}^n a_ix_i\Big|\geq \frac{1}{2}\sum_{i=1}^n a_iz_i\Big\}\Big|\geq
2^{-n}|K_p|\geq (2C_1)^{-n}\Pr(X\in C_1K_p)/g(0),
\]
therefore choosing $z$ such that $\sum_{i=1}^n a_iz_i=\|a\|_{K_p^o}$ we get
\begin{align*}
\Big\|\sum_{i=1}^n a_iX_i\Big\|_p&\geq 
2^{-1/p}\|a\|_{K_p^o} e^{-1}g(0)^{1/p}
\Big|\Big\{x\in K_p \colon \Big|\sum_{i=1}^n a_ix_i\Big|\geq \frac{1}{2}\sum_{i=1}^n a_iz_i\Big\}\Big|^{1/p}
\\
&\geq 2^{-1/p}\|a\|_{K_p^o} e^{-1}(2C_1)^{-n/p}\Pr(X\in C_1K_p)^{1/p}\geq
\frac{1}{4eC_1}\|a\|_{K_p^o}.
\end{align*}

To get the upper estimate  notice that
\[
\Pr\Big(\Big|\sum_{i=1}^n a_iX_i\Big|> C_1\|a\|_{K_p^o}\Big)\leq \Pr(X\notin C_1K_p)\leq e^{-p}.
\]
Together with the symmetry and log-concavity of $\sum_{i=1}^n a_iX_i$ this gives 
\[
\Pr\Big(\Big|\sum_{i=1}^n a_iX_i\Big|> C_1t\|a\|_{K_p^o}\Big)\leq  e^{-tp} \mbox{ for }t\geq 1.
\]
Integration by parts yields
$\Big\|\sum_{i=1}^n a_iX_i\Big\|_p\leq C\|a\|_{K_p^o}$.
\end{proof}

{\bf Remark.} The same argument as above shows that if $\alpha\geq e$ then
\[
\Big\|\sum_{i=1}^n a_iX_i\Big\|_p\geq 
\frac{1}{4\alpha C_1}\sup\Big\{\sum_{i=1}^n a_ix_i\colon g(x)\geq \alpha^{-p}g(0)\Big\}.
\]

From now on till the end of this section we  assume that vector $X$ is unconditional, log-concave and isotropic. Jensen's inequality
and Hitczenko  estimates for moments of Rademacher sums \cite{H} (see also \cite{MS}) imply that
for $p\geq 2$,
\begin{align}
\notag
\Big\|\sum_i a_iX_i\Big\|_p&=\Big\|\sum_i a_i\ve_i|X_i|\Big\|_p\geq \Big\|\sum_i a_i\ve_i\Ex|X_i|\Big\|_p
\\
\label{belowrad}
&\geq\frac{1}{C} \Big(\sum_{i\leq p}a_i^*+\sqrt{p}\Big(\sum_{i>p} |a_i^*|^2\Big)^{1/2}\Big).
\end{align}
The result of Bobkov and Nazarov \cite{BN} yields for $p\geq 2$,
\begin{align}
\label{upexp}
\Big\|\sum_i a_i X_i\Big\|_p\leq C 
\Big\|\sum_i a_i E_i\Big\|_p\leq C\Big(p\max_{i}|a_i|+\sqrt{p}\Big(\sum_{i} a_i^2\Big)^{1/2}\Big),
\end{align}
where $(E_i)$ is a sequence of independent symmetric exponential random variables with variance 1
and to get the second inequality we used the result of Gluskin and Kwapie\'n \cite{GK}.

Estimates \eqref{belowrad} and \eqref{upexp} together with Proposition \ref{largep} give
\begin{equation}
\label{inclKp}
\frac{1}{C}(\sqrt{p}B_2^n\cap B_{\infty}^n)\subset
\Big\{x\colon g(x)\geq e^{-p}g(0)\Big\}\subset C(\sqrt{p}B_2^n+pB_1^n)
\quad \mbox{ for }p\geq n.
\end{equation}

\begin{cor}
\label{unc_largep}
Let $X=(X_1,\ldots,X_n)$ be an unconditional log-concave isotropic vector  with the density $g$. 
Then for any $p\geq n$ we have
\begin{align*}
\Big\|\sum_{i=1}^n a_iX_i\Big\|_p
&\sim \sup\Big\{\sum_{i=1}^n a_ix_i\colon g(x)\geq e^{-p}g(0)\Big\}
\\
&\sim \sup\Big\{\sum_{i=1}^n a_ix_i\colon g(x)\geq e^{-5p/2}\Big\}
\\
&\sim \sup\Big\{\sum_{i=1}^n |a_i|t_i\colon \Pr(|X_1|\geq t_1,\ldots,|X_n|\geq t_n)\geq e^{-p}\Big\}.
\end{align*}
\end{cor}

\begin{proof}
We have $g(0)=L_{X}^n$, where $L_X$ is the isotropic constant of vector $X$. Unconditionality  
of $X$ implies boundedness of $L_X$, thus
\[
e^{-3n/2}\leq (2\pi e)^{-n/2}\leq g(0)\leq C_2^n,
\]
where $C_2$ is an absolute constant (see for example \cite{BN}). Hence
\begin{equation}
\label{incl}
\{x\colon g(x)\geq e^{-p}g(0)\} \subset \{x\colon g(x)\geq e^{-5p/2}\}\subset
\{x\colon g(x)\geq (e^{5/2}C_2)^{-p}g(0)\}
\end{equation}
and first two estimates on moments follows by Proposition \ref{largep} (see also remark after it).

For any $t_1,\ldots,t_n\geq 0$,
\[
\Ex\Big|\sum_{i=1}^na_iX_i\Big|^p\geq 
\Big(\sum_{i=1}^n|a_i|t_i\Big)^p2^{-n}\Pr(|X_i|\geq t_1,\ldots,|X_n|\geq t_n),
\]
therefore
\[
\Big\|\sum_{i=1}^n a_iX_i\Big\|_p\geq 
\frac{1}{2e}\sup\Big\{\sum_{i=1}^n a_it_i\colon \Pr(|X_1|\geq t_1,\ldots,|X_n|\geq t_n)\geq e^{-p}\Big\}.
\]
To prove the opposite estimate we use already proven bound and take $x$ such that $g(x)\geq e^{-5p/2}$ and 
$\sum_{i=1}^na_ix_i\geq \frac{1}{C_2}\|\sum_{i=1}^n a_iX_i\|_p$. By the unconditionality without loss of generality we may assume that all $a_i$'s and $x_i$'s are nonnegative. Notice that by \eqref{inclKp} and
\eqref{incl} we have $g(1/C_3,\ldots,1/C_3)\geq e^{-5p/2}$. Hence by log-concavity of $g$ we also
have $g(y)\geq e^{-5p/2}$ for $y_i=(x_i+1/C_3)/2$. Notice that $g$ is coordinate increasing
on $\er_+^n$, therefore
\[
\Pr\Big(X_1\geq \frac{y_1}{2},\ldots,X_n\geq \frac{y_n}{2}\Big)\geq g(y)\prod_{i=1}^n\frac{y_i}{2}
\geq e^{-5p/2}(4C_3)^{-n}\geq (4e^{5/2}C_3)^{-p}.
\] 
Function $F(s_1,\ldots,s_n)=-\ln\Pr(X_1\geq s_1,\ldots,X_n\geq s_n)$ is convex on $\er_+^n$, 
$F(0)=n\ln 2$, therefore
\[
\Pr\Big(|X_1|\geq \frac{y_1}{C_4},\ldots,|X_n|\geq \frac{y_n}{C_4}\Big)=
2^n\Pr\Big(X_1\geq \frac{y_1}{C_4},\ldots,X_n\geq \frac{y_n}{C_4}\Big)\geq e^{-p}
\]
for suffiently large $C_4$. To conclude it is enough to notice that
\[
\sum_{i=1}^n a_i\frac{y_i}{C_4}\geq \frac{1}{2C_4}\sum_{i=1}^n a_ix_i\geq 
\frac{1}{2C_2C_4}\Big\|\sum_{i=1}^n a_iX_i\Big\|_p.
\]
\end{proof}

\begin{thm}
\label{momunc}
Suppose that $X$ is an unconditional log-concave isotropic random vector in $\er^n$. 
Then for any $p\geq 2$,
\begin{align*}
\Big\|\sum_{i=1}^na_iX_i\Big\|_p&\sim
\sup\Big\{\sum_{i\in I_p} a_ix_i\colon\ g_{I_p}(x)\geq e^{-p}g_{I_p}(0)\Big\}+
\sqrt{p}\Big(\sum_{i\notin I_p}a_i^2\Big)^{1/2},
\\
&\sim \sup\Big\{\sum_{i\in I_p} a_ix_i\colon\ g_{I_p}(x)\geq e^{-5p/2}\Big\}+
\sqrt{p}\Big(\sum_{i\notin I_p}a_i^2\Big)^{1/2}
\\
&\sim \sup\Big\{\sum_{i\in I_p} |a_i|t_i\colon\ 
\Pr\Big(\forall_{i\in I_p}\ |X_i|\geq t_i\Big)\geq e^{-p}\Big\}+
\sqrt{p}\Big(\sum_{i\notin I_p}a_i^2\Big)^{1/2},
\end{align*}
where $g_{I_p}$ is the density of $(X_i)_{i\in I_p}$ and $I_p$ is the set of indices of
$\min\{\lceil p\rceil,n\}$ largest values of $|a_i|$'s.
\end{thm}

\begin{proof}
By Corollary \ref{unc_largep} it is enough to show that
\begin{align}
\notag
\frac{1}{C}\Big(\Big\|\sum_{i\in I_p}a_iX_i\Big\|_p+\sqrt{p}\Big(\sum_{i\notin I_p}a_i^2\Big)^{1/2}&\Big)
\leq \Big\|\sum_{i=1}^na_iX_i\Big\|_p
\\
\label{red_to_p}
&\leq
C\Big(\Big\|\sum_{i\in I_p}a_iX_i\Big\|_p+\sqrt{p}\Big(\sum_{i\notin I_p}a_i^2\Big)^{1/2}\Big).
\end{align}
Observe also that $\sum_{i\notin I_p}a_i^2=\sum_{i>p}|a_i^*|^2$.

Unconditionality of $X_i$ implies that $\|\sum_{i=1}^na_iX_i\|_p\geq \|\sum_{i\in I_p}a_iX_i\|_p$.
Hence the lower estimate in \eqref{red_to_p} follows by \eqref{belowrad}.

Obviously we have
\[
 \Big\|\sum_{i=1}^na_iX_i\Big\|_p\leq  \Big\|\sum_{i\in I_p}a_iX_i\Big\|_p+ 
\Big\|\sum_{i\notin I_p}a_iX_i\Big\|_p.
\]
Estimate \eqref{belowrad} and \eqref{upexp} imply
\begin{align*}
\Big\|\sum_{i\notin I_p}a_iX_i\Big\|_p&\leq 
C\Big(p\max_{i\notin I_p}|a_i|+\sqrt{p}\Big(\sum_{i\notin I_p}a_i^2\Big)^{1/2}\Big)
\\
&\leq C\Big(\Big\|\sum_{i\in I_p}a_iX_i\Big\|_p+\sqrt{p}\Big(\sum_{i\notin I_p}a_i^2\Big)^{1/2}\Big)
\end{align*}
and upper bound in \eqref{red_to_p} follows.
\end{proof}


{\bf Example 1.} Let $X_i$ be independent symmetric log-concave r.v's. Define 
$N_i(t):=-\Pr(|X_i|\geq t)$, then $\Pr(X_i\geq t_i \mbox{ for }i\in I_p)=\exp(-\sum_{i\in I_p}N_i(t_i))$
and Theorem \ref{momunc} yields the Gluskin-Kwapie\'n estimate
\[
\Big\|\sum_{i=1}^na_iX_i\Big\|_p\sim 
\sup\Big\{\sum_{i\in I_p} |a_i|t_i\colon \sum_{i\in I_p} N_i(t_i)\leq p\Big\}
+\sqrt{p}\Big(\sum_{i\notin I_p}a_i^{2}\Big)^{1/2}.
\]

\medskip

{\bf Example 2.} Let $X$ be uniformly distriputed on $r_{n,q}B_q^n$ with $1\leq q<\infty$, where $r_{n,q}$ is chosen in such a way that $X$ is isotropic. Then it is easy to check that $r_{n,q}\sim n^{1/q}$.
Since all $k$-dimensional sections of $B_q^n$ are homogenous we immediately obtain that for
$I\subset\{1,\ldots,n\}$ and $x\in \er^I$,
$g_{I}(x)/g_{I}(0)=(1-(\|x\|_{q}/r_{n,q})^q)^{n-|I|}$.
Hence for $1\leq p\leq n/2$  we get that
\[
\sup\Big\{\sum_{i\in I_p}a_ix_i\colon g_{I_p}(x)\geq e^{-p}g_{I_p}(0)\Big\}\sim
\sup\Big\{\sum_{i\in I_p}a_ix_i\colon \|x\|_{q}\leq p \Big\}.
\] 
Since for $p\geq n/2$, 
$\|\sum_{i=1}^na_iX_i\Big\|_p\sim \|\sum_{i=1}^na_iX_i\Big\|_{n/2}$, we
recover the result from \cite{BGMN} and show that for $p\geq 2$,
\[
\Big\|\sum_{i=1}^na_iX_i\Big\|_p\sim \min\{p,n\}^{1/q}\Big(\sum_{i\leq p}|a_i^*|^{q'}\Big)^{1/q'}+
\sqrt{p}\Big(\sum_{i>p}|a_i^*|^{2}\Big)^{1/2},
\]
where $1/q'+1/q=1$.

\medskip

{\bf Remark.} In the case of vector coefficients the following conjecture seems reasonable.
There exists a universal constant
$C$ such that for any isotropic unconditional log-concave vector $X=(X_1,\ldots,X_n)$ any vectors 
$v_1,\ldots,v_n$ in a normed space $(F,\|\cdot\|)$ and $p\geq 1$,
\[
\Big(\Ex\Big\|\sum_{i=1}^n v_iX_i\Big\|^p\Big)^{1/p}\sim
\Big(\Ex\Big\|\sum_{i=1}^n v_iX_i\Big\|+
\sup_{\|\varphi\|_{*}\leq 1}\Big(\Ex\Big|\sum_{i=1}^n \varphi(v_i)X_i\Big|^p\Big)^{1/p}\Big).
\]
The nontrivial part is the upper bound for $(\Ex\|\sum_{i=1}^n v_iX_i\|^p)^{1/p}$.
It is known that the above conjecture holds if the space $(F,\|\cdot\|)$ has nontrivial cotype --
see \cite{La2} for this and some related results.

\medskip

{\bf Remark.} Let $S=\sum_{i=1}^n a_iX_i$, where $X$ is as in Theorem \ref{momunc}. Then 
$\Pr(|S|\geq e\|S\|_p)\leq e^{-p}$ by the Chebyshev's inequality. Moreover $\|S\|_{2p}\leq C\|S\|_p$
for $p\geq 2$, hence Paley-Zygmund inequality yields $\Pr(|S|\geq \|S\|_p/C)\geq \min\{1/C,e^{-p}\}$.
This way Theorem \ref{momunc} may be also used to get two-sided estimates for tails of $S$.

\section{Gaussian approximation of moments}

Let $\gamma_p=\|{\cal N}(0,1)\|_p=2^{p/2}\Gamma(\frac{p+1}{2})/\sqrt{\pi}$.
In \cite{La} it was shown that for independent symmetric random variables $X_1,\ldots,X_n$ with log-concave
tails (notice that log-concave symmetric random variables have log-concave tails) and variance 1,
\begin{equation}
\label{indapr}
\Big|\Big\|\sum_{i=1}^n a_iX_i\Big\|_p-\gamma_p\|a\|_2\Big|\leq p\|a\|_{\infty}
\quad \mbox{ for }a\in \er^n,\ p\geq 3
\end{equation}
(see also \cite{Li} for $p\in [2,3)$).
The purpose of this section is to discuss  similar statements for general log-concave isotropic vectors $X$.

The lower estimate of moments is easy. In fact it holds for more general class of unconditional vectors with bounded fourth moment.

\begin{prop}
Suppose that $X$ is an isotropic unconditional $n$-dimensional vector with finite fourth moment. Then
for any nonzero $a\in \er^n$ and $p\geq 2$,
\begin{align*}
\Big\|\sum_{i=1}^n a_iX_i\Big\|_p&\geq 
\gamma_p\|a\|_2-\frac{p}{\sqrt{2}\|a\|_2}\Big(\sum_{i=1}^na_i^4\Ex X_i^4\Big)^{1/2}
\\
&\geq  \gamma_p\|a\|_2-\frac{p}{\sqrt{2}}\max_{i}(\Ex X_i^4)^{1/2}\|a\|_{\infty}.
\end{align*}
\end{prop}

\begin{proof} Let us fix $p\geq 2$.
By the homogenity we may and will assume that $\|a\|_2=1$.

Corollary 1 in \cite{La} gives 
\[
\Big\|\sum_{i=1}^n b_i\ve_i\Big\|_p\geq \gamma_p\Big(\sum_{i\geq \lceil p/2\rceil}|b_i^*|^2\Big)^{1/2}
\quad \mbox{ for } b\in \er^n,
\]
where $(b_i^*)$ denotes the nonicreasing rearrangement of $(|b_i|)_{i\leq n}$. Therefore
\begin{align*}
\Big\|&\sum_{i=1}^n a_iX_i\Big\|_p^p=\Ex\Big|\sum_{i=1}^n a_i\ve_iX_i\Big|^p\geq
\gamma_p^p\Ex\Big(\sum_{i=1}^na_i^2X_i^2-\max_{\# I<p/2}\sum_{i\in I}a_i^2 X_i^2\Big)^{p/2}
\\
&\geq \gamma_p^p\Big(\Ex\Big(\sum_{i=1}^na_i^2X_i^2-\max_{\# I<p/2}\sum_{i\in I}a_i^2 X_i^2\Big)\Big)^{p/2}
= \gamma_p^p\Big(1-\Ex\max_{\# I<p/2}\sum_{i\in I}a_i^2 X_i^2\Big)^{p/2}.
\end{align*}
We have
\begin{align*}
\Ex\max_{\# I<p/2}\sum_{i\in I}a_i^2 X_i^2&\leq 
\Ex\max_{\# I<p/2}\sqrt{\# I}\Big(\sum_{i\in I}a_i^4X_i^4\Big)^{1/2}
\leq \sqrt{\frac{p}{2}}\Ex\Big(\sum_{i=1}^na_i^4X_i^4\Big)^{1/2}
\\
&\leq \sqrt{\frac{p}{2}}\Big(\sum_{i=1}^na_i^4\Ex X_i^4\Big)^{1/2}.
\end{align*}
Since $\sqrt{1-x}\geq 1-x$ for $x\geq 0$ and $\gamma_p\leq \sqrt{p}$ the assertion easily follows.
\end{proof}

Since $\Ex Y^4\leq 6$ for symmetric log-concave random variables $Y$ we immediately get the following.

\begin{cor}
\label{lowerest}
Let $X$ be an isotropic unconditional $n$-dimensional log-concave vector. Then
for any $a\in \er^n\setminus\{0\}$ and $p\geq 2$,
\[
\Big\|\sum_{i=1}^n a_iX_i\Big\|_p\geq 
\gamma_p\|a\|_2-\frac{p}{\|a\|_2}\Big(3\sum_{i=1}^na_i^4\Big)^{1/2}\geq 
\gamma_p\|a\|_2-\sqrt{3}p\|a\|_{\infty}.
\]
\end{cor}

Now we turn our attention to the upper bound. Notice that for unconditional vectors $X$ and $p\geq 2$,
\begin{equation}
\label{upp1}
\Big\|\sum_{i=1}^n a_iX_i\Big\|_p=\Big\|\sum_{i=1}^n a_i\ve_iX_i\Big\|_p\leq
\gamma_p\Big\|\Big(\sum_{i=1}^n a_i^2X_i^2\Big)^{1/2}\Big\|_p,
\end{equation}
where the last inequality follows by the Khintchine inequality with optimal constant \cite{Ha}.
First we will bound moments of
$(\sum_{i=1}^n a_i^2X_i^2)^{1/2}$ using the result of Klartag \cite{K}.

\begin{prop}
For any isotropic unconditional $n$-dimensional log-concave vector $X$, $p\geq 2$ and nonzero
$a\in\er^n$ we have
\[
\Big\|\sum_{i=1}^n a_iX_i\Big\|_p-\gamma_p\|a\|_2
\leq Cp^{5/2}\frac{1}{\|a\|_2}\Big(\sum_{i=1}^n |a_i|^4\Big)^{1/2}
\leq Cp^{5/2}\|a\|_{\infty}. 
\]
\end{prop}

\begin{proof}
By homogenity we may assume that $\|a\|_2=1$.
We have 
\[
\Big\|\Big(\sum_{i=1}^n a_i^2X_i^2\Big)^{1/2}\Big\|_{p}
\leq 1+\Big\|\Big(\Big(\sum_{i=1}^n a_i^2X_i^2\Big)^{1/2}-1\Big)_{+}\|_{p}.
\]
Notice that
\[
\sum_{i=1}^n a_i^2(X_i^2-1)=
\Big(\Big(\sum_{i=1}^n a_i^2 X_i^2\Big)^{1/2}-1\Big)\Big(\Big(\sum_{i=1}^n a_i^2 X_i^2\Big)^{1/2}+1\Big),
\]
thus
\[
\Big\|\Big(\Big(\sum_{i=1}^n a_i^2X_i^2\Big)^{1/2}-1\Big)_{+}\|_{p}\leq 
\Big\|\sum_{i=1}^n a_i^2(X_i^2-1)\Big\|_p.
\]
Lemma 4 in \cite{K} gives
\[
\Big\|\sum_{i=1}^n a_i^2(X_i^2-1)\Big\|_2^2=\Var\Big(\sum_{i=1}^n a_i^2X_i^2\Big)
\leq \frac{8}{3}\sum_{i=1}^n a_i^4\Ex X_i^4\leq 16 \sum_{i=1}^n a_i^4.
\]
Comparison of moments of polynomials with respect to log-concave distributions \cite{NSV} implies
\[
\Big\|\sum_{i=1}^n a_i^2(X_i^2-1)\Big\|_p\leq (Cp)^2\Big\|\sum_{i=1}^n a_i^2(X_i^2-1)\Big\|_2
\leq  Cp^2\Big(\sum_{i=1}^n a_i^4\Big)^{1/2}.
\]
\end{proof}

We may improve $p^{5/2}$ term if we assume some concentration properties of vector $X$. We say that
vector $X$ satisfies  \emph{exponential concentration} with constant $\kappa$ if 
\[
\Pr(X\in A)\geq \frac{1}{2}\ \Rightarrow \ \Pr(X\in A+tB_2^n)\geq 1-e^{-t/\kappa}.
\] 

\begin{prop}
Let $X$ be an isotropic unconditional  vector that satisfies exponential concentration
with constant $\kappa$. Then for any $p\geq 2$ and $a\in \er^n$,
\[
\Big\|\sum_{i=1}^n a_iX_i\Big\|_p\leq \gamma_p\|a\|_2+ C\kappa p^{3/2}\|a\|_{\infty}.
\] 
\end{prop}

\begin{proof}
Notice that 
\[
\sup\Big\{\Big(\sum_{i=1}^n a_i^2y_i^2\Big)^{1/2}\colon y\in tB_2^n\Big\}=t\|a\|_{\infty}.
\]
Using standard arguments we may therefore show that exponential concentration implies for $p\geq 2$,
\[
\Big\|\Big(\sum_{i=1}^n a_i^2X_i^2\Big)^{1/2}\Big\|_p\leq 
\Big\|\Big(\sum_{i=1}^n a_i^2X_i^2\Big)^{1/2}\Big\|_2+C\kappa p\|a\|_{\infty}
=\|a\|_2+C\kappa p\|a\|_{\infty}.
\]
We conclude using \eqref{upp1}.
\end{proof}

Since by the result of Klartag \cite{K} unconditional log-concave vectors satisfy exponential concentration with constant $C\log n$  we get

\begin{cor}
Let $X$ be isotropic unconditional  logconcave vector. Then for any $p\geq 2$ and $a\in \er^n$,
\[
\Big\|\sum_{i=1}^n a_iX_i\Big\|_p\leq \gamma_p\|a\|_2+ C p^{3/2}\log n\|a\|_{\infty}.
\] 
\end{cor}

To get the factor $p$ instead of $p^{3/2}$ we need a stronger notion than exponential concentration.
We say that a random vector $X$ satisfies \emph{two level concentration with constant $\kappa$} if
\[
\Pr(X\in A)\geq \frac{1}{2}\ \Rightarrow \ \Pr(X\in A+\sqrt{t}B_2^n+tB_1^n)\geq 1-e^{-t/\kappa}.
\] 

Since it is enough to consider $t\geq 1$ two level concentration is indeed stronger than exponential
concentration. 

\begin{prop}
Suppose that $X$ is an isotropic unconditional vector that satisfies two level concentration
with constant $\kappa$. Then for any $p\geq 2$ and $a\in \er^n$,
\[
\Big\|\sum_{i=1}^n a_iX_i\Big\|_p\leq \gamma_p\|a\|_2+ C\kappa p\|a\|_{\infty}.
\] 
\end{prop}

\begin{proof}
For $p\geq 2$ define a norm $\trip\cdot\trip_p$ on $\er^n$ by 
$\trip x\trip_p=\|\sum_{i=1}^nx_i\ve_i\|_p$. Notice that $\trip x\trip_p\leq \gamma_p\|x\|_2$, hence
\[
\Ex\trip(a_iX_i)\trip_p^2\leq \gamma_p^2\|a\|_2^2. 
\]

Observe also that
\begin{align*}
\sup\{\trip(a_ix_i)\trip_p&\colon x\in \sqrt{t}B_2^n+tB_1^n\}
\\
&\leq
\sqrt{t}\sup\{\trip(a_ix_i)\trip_p\colon x\in B_2^n\}+t\sup_{j\leq n}\trip(a_{i}\delta_{i,j})\trip_p
\\
&\leq \sqrt{t}\gamma_p\sup\{\|(a_ix_i)\|_2\colon x\in B_2^n\}+t\|a\|_{\infty}=
(\sqrt{t}\gamma_p+t)\|a\|_{\infty}.
\end{align*}

Let $M_p=\Med(\trip(a_iX_i)\trip_p)$, two level concentration (applied twice to sets
$A=\{\trip(a_iX_i)\trip_p\leq M_p\}$ and $A=\{\trip(a_iX_i)\trip_p\geq M_p\}$) implies that
\[
\Pr\Big(\Big|\trip(a_iX_i)\trip_p-M_p\Big|\geq(\sqrt{t}\gamma_p+t)\|a\|_{\infty}\Big)\leq 2\exp(-t/\kappa).
\]
Integrating by parts this gives for $p\geq q\geq 2$,
\[
\big\|\trip(a_iX_i)\trip_p-M_p\big\|_q\leq C\kappa(\sqrt{q}\gamma_p+q)\|a\|_{\infty}
\leq C\kappa p\|a\|_{\infty}.
\]
Hence
\[
\Big\|\sum_{i=1}^n a_iX_i\Big\|_p=\|\trip(a_iX_i)\trip_p\|_p
\leq \|\trip(a_iX_i)\trip_p\|_2+C\kappa p\|a\|_{\infty}
\leq \gamma_p\|a\|_2+C\kappa p\|a\|_{\infty}.
\]
\end{proof}
 
Unfortunately we do not know many examples of random vectors satisfying two level concentration with
a good constant. Using estimate \eqref{upexp}
it is not hard to see that infimum convolution inequality investigated in \cite{LW} implies two level
concentration. In particular isotropic
log-concave unconditional vectors with independent coordinates and isotropic vectors uniformly
distributed on the (suitably rescaled) $B_p^n$ balls satisfy two level concentration with an absolute constant. 

The last approach to the problem of Gaussian approximation of moments we will discuss is based on 
the notion of negative association. We say that random variables
$(Y_1,\ldots,Y_n)$ are \emph{negatively associated} if for any disjoint sets $I_1, I_2$ in $\{1,\ldots,n\}$
and any bounded functions $f_i\colon \er^{I_i}\rightarrow \er$, $i=1,2$ that are coordinate nondecreasing
we have 
\[
\Cov\Big(f_1((X_i)_{i\in I_1}),f_2((X_i)_{i\in I_2})\Big)\leq 0.
\]

Our next result is an unconditional version of Theorem 1 in \cite{Sh}.

\begin{thm}
\label{negmod}
Suppose that $X=(X_1,\ldots,X_n)$ is an unconditional random vector with finite second moment and 
random variables $(|X_i|)_{i=1}^n$ are negatively associated. Let $X_1^*,\ldots,X_n^*$ be independent
random variables such that $X_i^*$ has the same distribution as $X_i$. Then for any nonnegative
function $f$ on $\er$ such that $f''$ is convex and any $a_1,\ldots,a_n$ we have
\begin{equation}
\label{negcomp}
\Ex f\Big(\sum_{i=1}^n a_iX_i\Big)\leq \Ex f\Big(\sum_{i=1}^n a_iX_i^*\Big).
\end{equation}
In particular
\[
\Ex \Big|\sum_{i=1}^n a_iX_i\Big|^p\leq \Ex \Big|\sum_{i=1}^n a_iX_i^*\Big|^p \quad \mbox{ for }p\geq 3.
\] 
\end{thm}

\begin{proof}
Since random variables $|a_iX_i|$ are also negatively associated, it is enough to consider the case when
$a_i=1$ for all $i$.
We may also assume that variables $X_i^*$ are independent of $X$. Assume first that random variables $X_i$
are bounded.


Let $Y=(Y_1,\ldots,Y_n)$ be independent copy of $X$ and $2\leq k\leq n$. To shorten the notation put for
$1\leq l\leq n$,
$S_l=\sum_{i=1}^l\ve_i|X_i|$ and $\tilde{S}_l=\sum_{i=1}^l\ve_i|Y_i|$ (recall that
$\ve_i$ denotes a Bernoulli sequence independent of other variables). 

 We have
\begin{align}
\notag
f(S_{k})&+f(\tilde{S}_k)
-f(S_{k-1}+\ve_k |Y_k|)-f(\tilde{S}_{k-1}+\ve_{k}|X_k|)
\\
\notag
&=\int_{|Y_k|}^{|X_k|}\ve_k (f'(S_{k-1}+\ve_k t)-f'(\tilde{S}_{k-1}+\ve_k t))dt
\\
\label{neg1}
&=\int_{-\infty}^{\infty}\ve_k (f'(S_{k-1}+\ve_k t)-f'(\tilde{S}_{k-1}+\ve_k t))
(\ind_{\{|X_k|\geq t\}}-\ind_{\{|Y_k|\geq t\}})dt.
\end{align} 

Define for $t>0$, $g_t(x)=\Ex \ve_k f'(x+\ve_k t)=(f'(x+t)-f'(x-t))/2$ and
\[
h_t(|x_1|,\ldots,|x_{k-1}|)=\Ex_{\ve}\ve_kf'\Big(\sum_{i=1}^{k-1} \ve_i|x_i|+\ve_k t\Big)=
\Ex g_t\Big(\sum_{i=1}^{k-1} \ve_i|x_i|\Big).
\]

Taking the expectation in \eqref{neg1} and using the unconditionality we get
\begin{align*}
2\Big(&\Ex f\Big(\sum_{i=1}^k X_i\Big)-\Ex f\Big(\sum_{i=1}^{k-1}X_i+X_k^*\Big)\Big)
\\
&=\Ex\int_{-\infty}^{\infty}\ve_k (f'(S_{k-1}+\ve_k t)-f'(\tilde{S}_{k-1}+\ve_k t))
(\ind_{\{|X_k|\geq t\}}-\ind_{\{|Y_k|\geq t\}})dt
\\
&=\int_{-\infty}^{\infty}\Ex \big[\big(h_t(|X_1|,\ldots,|X_k|)-h_t(|Y_1|,\ldots,|Y_k|)\big)
\big(\ind_{\{|X_n|\geq t\}}-\ind_{\{|Y_n|\geq t\}}\big)\big]dt
\\
&=\int_{-\infty}^{\infty}\Cov\big(h_t(|X_1|,\ldots,|X_{k-1}|),\ind_{\{|X_k|\geq t\}}\big)dt.
\end{align*}

Convexity of $f''$ implies that the function $g_t$ is convex on $\er$, therefore the 
function $h_t$ is coordinate increasing on $\er_+^{k-1}$. 
So by the negative association we get
\begin{equation}
\label{onestar}
\Ex f\Big(\sum_{i=1}^k X_i\Big)\leq \Ex f\Big(\sum_{i=1}^{k-1}X_i+X_k^*\Big)
\end{equation}
The same inequality holds if we change the function $f$ into the function $f(\cdot+h)$ for any
$h\in \er$. Therefore applying \eqref{onestar} conditionally we get
\[
\Ex f\Big(\sum_{i=1}^k X_i+\sum_{i=k+1}^n X_i^*\Big)
\leq \Ex f\Big(\sum_{i=1}^{k-1}X_i+\sum_{i=k}^nX_i^*\Big)
\]
and inequality \eqref{negcomp} easily follows in the bounded case.

To settle the unbounded case first notice that random variables $|X_i|\wedge m$ are bounded and negatively
associated for any $m>0$. Hence we know that
\[
\Ex f\Big(\sum_{i=1}^n \ve_i|X_i|\wedge m\Big)\leq \Ex f\Big(\sum_{i=1}^{n}\ve_i|X_i^*|\wedge m\Big).
\] 
We have $\liminf_{m\rightarrow\infty}f(\sum_{i=1}^n \ve_i|X_i|\wedge m)\geq \Ex f(\sum_{i=1}^n \ve_i|X_i|)$,
so it is enough to show that , 
$\liminf_{m\rightarrow\infty}\Ex f(\sum_{i=1}^n \ve_i|X_i^*|\wedge m)\leq \Ex f(\sum_{i=1}^n \ve_i|X_i^*|)$.

Let us define  $u(x)=f(x)-\frac{1}{2}f''(0)x^2$, function $u''$ is convex and $u''(0)=0$.
 Since $\Ex|X_i|^2=\Ex|X_i^*|^2<\infty$ it is enough to show that
for any $m>0$,
\begin{equation}
\label{aim}
\Ex u\Big(\sum_{i=1}^{n}\ve_i|X_i^*|\wedge m\Big)\leq \Ex u\Big(\sum_{i=1}^{n}\ve_i|X_i^*|\Big).
\end{equation}
Let for $s\in \er$, $v_s(t):= \Ex u(\ve_1 s+\ve_2 t)$, then $v_s''(t)=\Ex u''(\ve_1 s+\ve_2 t)\geq u''(\Ex(\ve_1 s+\ve_2 t))=0$ and $v_s'(0)=0$, hence $v_s$ is nondecreasing on $[0,\infty)$. 
Thus for any $x\in \er^n$,
\[
\Ex_{\ve}u\Big(\sum_{i=1}^{n}\ve_i|x_i|\wedge m\Big)\leq \Ex_{\ve} u\Big(\sum_{i=1}^{n}\ve_i|x_i|\Big)
\]
and \eqref{aim} immediately follows.
\end{proof}

\begin{cor}
Suppose that $X$ is an isotropic undonditional $n$-dimensional log-concave vector such that
variables $|X_i|$ are negatively associated. Then for any $a_1,\ldots,a_n$ and $p\geq 3$,
\[
-\sqrt{3}p\|a\|_{\infty}\leq \Big\|\sum_{i=1}^n a_iX_i\Big\|_p-\gamma_p\|a\|_2\leq p\|a\|_{\infty}.
\]
In particular the above inequality holds if $X$ has a uniform distribution on a (suitably rescaled)
Orlicz ball.
\end{cor}

\begin{proof}
First inequality follows by Corollary \ref{lowerest}, second by Theorem \ref{negmod} and \eqref{indapr}.
The last part of the statement is a consequence of the result of Pilipczuk and Wojtaszczyk \cite{PW}
(see also \cite{W} for a simpler proof and a slightly more general class of unconditional log-concave
measures with negatively associated absolute values of coordinates).
\end{proof}




{\bf Acknowledgements.} Part of the work was done at the Newton institute for Mathematical Sciences in Cambridge (UK) during the program ”Discrete Analysis”.

\end{document}